\documentstyle[12pt]{article}
\newcommand{\const}{\mathop{\rm const}\limits}

\begin{document}

\begin{center}

\vspace{3mm}

{\bf INDIVIDUAL LOWER BOUND FOR CALDERON'S }\par

\vspace{4mm}

{\bf GENERALIZED LORENTZ NORM ESTIMATES.} \\

\vspace{4mm}

 $ {\bf E.Ostrovsky^a, \ \ L.Sirota^b } $ \\

\vspace{4mm}

$ ^a $ Corresponding Author. Department of Mathematics and computer science, Bar-Ilan University, 84105, Ramat Gan, Israel.\\
\end{center}
E - mail: \ galo@list.ru \  eugostrovsky@list.ru\\
\begin{center}
$ ^b $  Department of Mathematics and computer science. Bar-Ilan University,
84105, Ramat Gan, Israel.\\
\end{center}
E - mail: \ sirota3@bezeqint.net\\

\vspace{3mm}
\begin{center}
                    {\sc Abstract.}\\

 \end{center}

 \vspace{4mm}

 We find the exact values for constants in bilateral Calderon - Stein - Weiss inequalities between  tail (Marcinkiewicz) norm
and weak Lebesgue (Lorentz) norm.\par

 Possible applications: Functional Analysis (for instance, interpolation of operators), Integral Equations,
Probability Theory and Statistics (tail estimations for random variables) etc.\par

 \vspace{4mm}

 \par

\vspace{3mm}


\vspace{3mm}

\section{Notations. Statement of problem.}

\vspace{3mm}

 Let $ (X = \{x\}, \cal{A}, \mu) $ be measurable space with  non-trivial sigma-finite measure $ \mu. $
We will suppose without loss of generality in the case  $ \mu(X) < \infty $  that $ \mu(X) = 1 $ (the probabilistic case)
and denote $  x = \omega, \ {\bf P} = \mu. $\par
 Define as usually for arbitrary measurable function $ f: X \to R $ its distribution function
 (more exactly, tail function)

 $$
 T_f(t) = \mu\{x: |f(x)| \ge t \}, \  t \ge 0,
 $$

 $$
 ||f||_p = \left[ \int_X |f(x)|^p \ \mu(dx) \right]^{1/p}, \ p \ge 1; \ L_p = \{f, ||f||_p < \infty\},
 $$
 and denote by $ f^*(t) = T_f^{-1}(t) $ the right inverse to the tail function $ T_f(t): $

 $$
 f^*(t) = \inf \{s: \ \mu( \{x: |f(x)| > s \} ) \le t \}.
 $$
   The following function $ f^{**}= f^{**}(t) $ play a very important role in the theory of interpolation of
 operators and harmonic analysis, see \cite{Bennet1}, \cite{Stein1}:
$$
 f^{**}(t) \stackrel{def}{=} t^{-1} \int_0^t f^*(s) \ ds, \ t > 0.
$$

  We will denote the set of all tail functions  as $ \{ T \}; $ obviously, the set  $ \{ T \}$ contains on
all the functions $ \{ H = H(t), \ t \ge 0 \}  $ which are right continuous, monotonically non-increasing
with values in the set  $ [0, \mu(X)]. $ \par

 Let  $ w = w(s), s \ge 0 $ be any (measurable) non-increasing numerical function (weight) defined on the set
 $ s \in (0, \infty) $  such that
 $$
 w(s) = 0 \Leftrightarrow s = 0;  \ \lim_{s \to \infty}w(s) = \infty.\eqno(1.1)
 $$

 The set of all such a functions we will denote as $ V; \ V = \{w\}. $  \par

  Moreover, we introduce the set  of all a weight functions $ W = \{ w \} $ under another  restriction:

$$
\forall w \in W \ \exists T \in \{ T \}  \Rightarrow w(T(s)) = 1/s. \eqno(1.2)
$$

 Let us introduce the following important functional

 $$
 \gamma(w) = \sup_{t > 0} \left[  \frac{w(t)}{t} \ \int_0^t \frac{du}{w(u)} \right] \eqno(1.3)
 $$
and the following quasi-norms:

$$
||f||_w^* = \sup_{t > 0} [w(t) \ f^*(t)], \eqno(1.4)
$$

$$
||f||_w = \sup_{t > 0} [w(t) \ f^{**}(t)], \eqno(1.5)
$$
 The necessary and sufficient condition for finiteness of the functional $ \gamma(w) $ and following for the
normability of the space $ L_w $  see, e.g. in the article \cite{Astashkin1}; see also \cite{Cwikel3}, \cite{Kaminska1}.  \par

\vspace{3mm}

{\bf Remark 1.1.} As long as
\vspace{3mm}
$$
f^{**}(t)= t^{-1} \sup_{\mu(E) \le t} \int_E |f(x)| \ \mu(dx),
$$
we can rewrite  the expression for $ ||f||_w  $ as follows:

$$
||f||_w = \sup_{t > 0} \left[ (w(t)/t) \cdot \sup_{E: \mu(E) \le t} \int_E  |f(x)| \ \mu(dx)   \right]. \eqno(1.6)
$$

 If the measure $ \mu $ has not atoms, then the expression (1.6) may be rewritten as follows:

 $$
||f||_w = \sup_{ E: 0 < \mu(E) < \infty } \left[ \frac{ w(\mu(E))}{\mu(E)} \cdot \int_E |f(x)| \ \mu(dx) \right]. \eqno(1.7)
 $$

  It follows from  identity (1.6) that $ ||f||_w $ is  rearrangement invariant norm and the space
$ L_w = \{ f: ||f||_w < \infty  \}  $ is (complete) Banach functional rearrangement invariant
 space with Fatou property.  The proof is alike to one in the case $  w(t) = t^{1/p}, \ p \ge 1;  $ see \cite{Bennet1}, chapters 1,2;
\cite{Stein1}, chapter 5, section 3. \par

\vspace{3mm}

 The norm $  ||f||_w   $ is named Marcinkiewicz's norm, see \cite{Krein1},  chapter 2, section 2.  More information about
considered in this article  Marcinkiewicz's and Lorentz (weak Lebesgue) spaces  with described applications see, e.g. in
\cite{Astashkin1}, \cite{Bennet1}, chapter 3, section 3; \cite{Calderon1},
\cite{Carro1}, \cite{Cwikel2}, \cite{Krein1}, \cite{Kufner1}, \cite{Lieb1}, \cite{Lorentz1}, \cite{Okikiolu1}, chapter 5,  section 5;
\cite{Soria1}, \cite{Stein1}, chapter 5, section 3 etc. \par

 See also many works of L.Maligranda  at all \cite{Cwikel3}, \cite{Kaminska1}, \cite{Kufner2} etc; M.M.Milman at all
\cite{Milman1}, \cite{Milman2}, \cite{Milman3}  etc. which are devoted to the theory of those spaces. \par

\vspace{3mm}

\section{Main result.}

\vspace{3mm}

 In the article \cite{Ostrovsky9} (Theorem 2.1.) has been proved the following estimation, which is some generalization of
 Calderon - Stein - Weiss bilateral inequality; see also \cite{Calderon1} and
 the  classical monographs  of E.M. Stein - G.Weiss \cite{Stein1}, chapter 5, section 3 and G.O.Okikiolu
 \cite{Okikiolu1},  chapter 5, section 5;  if

$$
w \in W, \ \gamma(w) < \infty,  \eqno(2.1)
$$
then

$$
1 \cdot||f||^*_w  \le ||f||_w \le \gamma(w) \cdot ||f||^*_w, \eqno(2.2)
$$
 and both  the coefficients $ "1" $ and $ "\gamma(w)" $ in (2.2) are the best possible. \par
\vspace{3mm}
 {\bf Remark 2.1.} Note that more general version of used in the proving of assertion (2.2)  Hardy's  inequality
 is obtained, for example, in  \cite{Arino1}. \par
 \vspace{3mm}
 The exactness of the constant $ " 1 " $ follows immediately from the consideration of the case $ w = w(s)= w_p(s) := s^{1/p}, \ p > 1. $
Namely, in this case we obtain the classical inequality

$$
||f||^*_{w_p}  \le ||f||_{w_p} \le \frac{p}{p-1} \cdot ||f||^*_{w_p},
$$
see  \cite{Stein1}, chapter 5, section 3; note that  $ \lim_{p \to \infty} p/(p-1)=1.  $\par

\vspace{3mm}

{\it  Notice that the exactness of lower bound in (2.2) is understood over all the set } $  W, $ \\
{\it in contradiction  to the upper estimation,  which is true for every function } $  w, \ w \in W. $ \par

\vspace{3mm}

{\bf  The aim of this short report is to prove that the lower coefficient "1" in the left hand-side of
bilateral inequality (2.2) is also exact (under some natural conditions)   for each function $  w, \ w \in V. $ }\par

\vspace{3mm}

 It suffices to consider only the probabilistic case  $ \mu(X)  = {\bf P}(X) = 1. $ \par
\vspace{3mm}

 More detail, let us denote

$$
\Theta(w) = \inf_{f \ne 0} \left[ \frac{||f||_w}{||f||^*_w} \right], \eqno(2.3)
$$

$$
G_{\kappa}(w) = \sup_{t \in (0,1)} \left[ w(t)(1-t^{\kappa} ) \right], \ \kappa = \const \in (0, \infty), \eqno(2.4)
$$

$$
H_{\kappa}(w) = \sup_{t \in (0,1)} \left[ w(t)(1-t^{\kappa}/(\kappa+1) ) \right]. \eqno(2.5)
$$

$$
K_{\kappa} (w)= \frac{G_{\kappa}(w)}{H_{\kappa}(w)},  \    \ K(w) =  \inf_{\kappa \in (0,\infty)}K_{\kappa} (w). \eqno(2.6)
$$

\vspace{3mm}

{\bf Theorem 1.}

$$
\Theta(w) \le K(w). \eqno(2.7)
$$

\vspace{3mm}

{\bf Proof.} Let $ \kappa $ be a fix number from the semi-axes $ (0,\infty). $  There exists a measurable function
$ f_{\kappa} = f_{\kappa}(x), \ x \in X $
(random variable)  which  may be defined for instance on the set $ X = \{x\}= [0,1] $ equipped with Lebesgue measure $ {\bf P} $ such that

$$
f^*_{\kappa}(t) = 1-t^{\kappa}, \ t \in [0,1], \eqno(2.8)
$$
for example
$$
f_{\kappa}(x):= 1 - x^{\kappa}, \ x \in [0,1];
$$
then

$$
f^{**}_{\kappa}(t) = 1-t^{\kappa}/ (\kappa+1), \ t \in [0,1]. \eqno(2.9)
$$

 Note that for all the values $  \kappa > 0 $

$$
\Theta(w)  \le \frac{||f_{\kappa}||_w}{||f_{\kappa}||^*_w}
=  \frac{\sup_{t \in (0,1)} \left[ w(t)(1-t^{\kappa} ) \right]}{\sup_{t \in (0,1)} \left[ w(t)(1-t^{\kappa} /( \kappa+1) ) \right]}=
\frac{G_{\kappa}(w)}{H_{\kappa}(w)} = K_{\kappa}(w), \eqno(2.10)
$$
therefore

$$
\Theta(w) \le \inf_{\kappa \in (0,\infty)} K_{\kappa}(w) = K(w). \eqno(2.11)
$$

\vspace{3mm}

\section{Examples.}

\vspace{3mm}
 We assume in this section $ w \in V. $\par

 {\bf Proposition 1.} Let $ w = w(s)= w_p(s) = s^{1/p}, \ p = \const > 1; $ then
 $$
 \Theta(w_p) = 1.\eqno(3.1)
 $$

{\bf Proof.} The lower bound $ \Theta(w) \ge 1, \ w \in V $ is obvious, see e.g. \cite{Stein1}, chapter 5, section 3, as long as
$ f^{**}(t) \ge f^*(t).  $ \par

 In order to prove the
upper estimate $ \Theta(w_p) \le 1 $ we use the assertion of theorem 1. Namely, let $  \kappa = \const \in (0,\infty); $ note that

$$
 \lim_{\kappa \to \infty} \left[  \frac{ f^*_{\kappa}(t)}{ f^{**}_{\kappa}(t)} \right] = 1, \ t \in (0,1);
$$
hence it is naturally to hope that

$$
 \lim_{\kappa \to \infty}  \frac{G_{\kappa}(w_p)}{H_{\kappa}(w_p)} = 1. \eqno(3.2)
$$

 In detail, we find by direct computations:

 $$
 G_{\kappa}(w_p) = \sup_{t \in (0,1)}  \left( t^{1/p} - t^{\kappa + 1/p }  \right) =
 (\kappa p + 1)^{-1/(\kappa p)} \cdot \frac{\kappa p}{\kappa p + 1}; \eqno(3.3)
 $$

 $$
 H_{\kappa}(w_p) = \sup_{t \in (0,1)}  \left( t^{1/p} - t^{\kappa + 1/p }/(\kappa+1)  \right) =
 \frac{ (\kappa + 1)^{1/(\kappa p)}}{(\kappa p + 1)^{1/(\kappa p)}} \cdot \frac{\kappa p}{\kappa p + 1}; \eqno(3.4)
 $$

$$
\frac{G_{\kappa}(w_p)}{H_{\kappa}(w_p)} = (\kappa + 1)^{1/(\kappa p)}, \eqno(3.5)
$$

$$
 K(w) =  \inf_{\kappa \in (0,\infty)}K_{\kappa} (w) \le \lim_{\kappa \to \infty}
K_{\kappa} (w)= \lim_{\kappa \to \infty}\frac{G_{\kappa}(w)}{H_{\kappa}(w)} =
$$

$$
\lim_{\kappa \to \infty}(\kappa + 1)^{1/(\kappa p)} = 1.     \eqno(3.6)
$$

\vspace{3mm}
{\bf Remark 3.1.} At the same result as in proposition 1 is true for the functions of a view

$$
w_{p,q}(s) = s^{1/p} \ (|\log s| + 1)^{1/q}, \ p > 1, q \in R; \eqno(3.7)
$$

\vspace{3mm}

$$
w_{p,q,r}(s) = s^{1/p} \ (|\log s| + 1)^{q}  \ ( \log(|\log s| + 3  )^{r}, \ p > 1, q,r \in R \eqno(3.8)
$$
etc.\par

\end{document}